\documentclass[a4paper]{amsart}
\usepackage{amsmath,amssymb}

\theoremstyle{definition}

\newcommand{\N}{\mathbb N}

\newtheorem*{propertyd0}{Property D0}
\newtheorem*{propertyd}{Property D}

 \DeclareMathOperator{\supp}{supp}
 
\DeclareMathOperator{\Ker}{Ker}

\newtheorem{thm}{Theorem}[section]

\newtheorem{prop}[thm]{Proposition}

\newtheorem{conj}[thm]{Conjecture}
\newtheorem{lem}[thm]{Lemma}

\begin{document}



\title[On the Erd{\H o}s--Ginzburg--Ziv constant]{On the Erd{\H o}s--Ginzburg--Ziv constant \\ of finite
abelian groups of high rank}

\address{Center for Combinatorics, LPMC-TJKLC,
Nankai University, Tianjin 300071, P.R. China}
\email{fys850820@163.com, wdgao1963@yahoo.com.cn,
zhongqinghai@yahoo.com.cn}

\author{Yushuang Fan and Weidong Gao  and Qinghai Zhong}

\maketitle

\begin{abstract}
Let $G$ be a finite abelian group. The Erd{\H o}s--Ginzburg--Ziv constant $\mathsf s (G)$ of $G$ is defined as the smallest integer $l \in \mathbb N$ such that every sequence \ $S$ \ over $G$ of length $|S| \ge l$ \ has a  zero-sum subsequence $T$ of length $|T| = \exp (G)$. If $G$ has rank at most two, then the precise value of $\mathsf s (G)$ is known (for cyclic groups this is the  Theorem of Erd{\H o}s--Ginzburg--Ziv). Only very little is known for groups of higher rank. In the present paper, we focus on
groups of the form  $G = C_n^r$,  with $n, r \in \N$ and $n \ge 2$,
and we tackle the study of $\mathsf s (G)$ with a new approach,
combining the direct problem with the associated inverse problem.
\end{abstract}

\bigskip
\section{Introduction and Main Result}
\bigskip

Let $G$ be an  additive finite abelian group. We denote by
\begin{itemize}
\item \ $\mathsf D (G)$ \  the smallest integer $l \in \mathbb
      N$ such that every sequence \ $S$ \ over $G$ of length $|S| \ge l$ \ has a non-empty zero-sum
      subsequence.

\smallskip
\item \ $\mathsf s (G)$ \  the smallest integer $l \in \mathbb
      N$ such that every sequence \ $S$ \ over $G$ of length $|S| \ge l$ \ has a  zero-sum
      subsequence $T$ of length $|T| = \exp (G)$.
\end{itemize}
Then $\mathsf D (G)$ is called the {\it Davenport constant} and
$\mathsf s (G)$ the {\it Erd{\H o}s-Ginzburg-Ziv constant} of $G$.
These are classical invariants in Combinatorial Number Theory, and
their precise values are known for groups with rank at most two.
Indeed, we have (see \cite[Theorem 5.8.3]{Ge-HK06a})

\medskip
{\bf Theorem A.} {\it  Let  $G = C_{n_1} \oplus C_{n_2}$ with $1 \le
n_1 \mid n_2$. Then
\[
\mathsf D (G) = n_1+n_2-1 \quad \text{and} \quad \mathsf s (G) =
2n_1+2n_2-3 \,.
\]
}

The result for $\mathsf D (G)$ dates back to the 1960s, and the
special case $n_1=1$ and $\mathsf s (C_{n_2}) = 2n_2-1$ is the
well-known Theorem of Erd{\H o}s-Ginzburg-Ziv proved in 1961
(\cite{Er-Gi-Zi61}). However, the special case where $n_1=n_2$ is a
prime was only settled in 2007 by C. Reiher (\cite{Re07a}). More
information can be found in the surveys \cite{Ga-Ge06b, Ge09a}.
Both the Davenport constant and the Erd{\H o}s-Ginzburg-Ziv
constant have found far reaching generalizations, and for  these
generalized versions, the precise values have been determined for
groups with rank at most two (see \cite[Section 6.1]{Ge-HK06a},
\cite{Fr-Sc10a}, \cite[Theorem 5.2]{Ge-Gr-Sc11c}, \cite{Pl-Sc11a}).

\smallskip
The situation is very different for groups of higher rank. Even for
the group $G = C_n \oplus C_n \oplus C_n$ with $n \ge 2$, the
precise value of the Davenport constant is unknown (for general $n$)
and the same is true for the Erd{\H o}s-Ginzburg-Ziv constant. In
what follows, we focus our discussion on the Erd{\H o}s-Ginzburg-Ziv
constant, which will be the main topic of the present paper. In 1995,
N. Alon and M. Dubiner \cite{Al-Du93} proved  that for every
positive integer $r$ there is a constant $c(r)$ depending only on
$r$ such that $\mathsf s(C_n^r)\leq c(r)n$ for all $n \ge 2$. To
illustrate the difficulties for obtaining precise values, let us
consider the special case $G = \mathbb F_3^r$, where $\mathbb
F_3$ is the finite field with three elements. Then $(\mathsf s
(G)-1)/2$ equals the maximal size of a cap in the affine space
$\mathbb F_3^r$. The maximal size of such caps has been studied in
finite geometry for decades, and the precise value is known so far
only for $r \le 6$ (see \cite{Po08a, Ed08a}). The connection to
affine caps will be addressed in greater detail in Section \ref{4}.
In the next theorem, we  gather the cases where precise values for
$\mathsf s (G)$ are known (more on upper and lower bounds will be
given in Section \ref{2}).

\medskip
{\bf Theorem B.} {\it Let $G$ be a finite abelian group, $n,r$
positive integers,   and  $a,b$ nonnegative integers.

\begin{itemize}
\smallskip
\item If $G = C_{2^a} \oplus C_{2^b}^{r-1}$ where $r \ge 2$, $b \ge 1$ and $a
\in [1, b]$, then $\mathsf s (G) = 2^{r-1}(2^a + 2^b - 2) + 1$
$($\cite[Corollary 4.4]{E-E-G-K-R07}$)$.

\smallskip
\item $\mathsf s(C_{3^a5^b}^3)=9(3^a5^b-1)+1$, where $a+b\geq 1$
$($\cite[Theorem 1.7]{G-H-S-T07}$)$.

\smallskip
\item $\mathsf s(C_{3^a}^4)=20(3^a-1)+1$, where $a\geq 1.$ {\rm (the precise value for $\mathsf s (C_3^4)$ was found independently  several times, see
  \cite[Section 5]{E-E-G-K-R07}; then use \cite[Theorems 1.3 and
  1.4]{E-E-G-K-R07})}.

\smallskip
\item $\mathsf s(C_3^5)=91$ and $\mathsf s(C_3^6)=225$
{\rm (see \cite[Theorem 1.2]{E-F-L-S02}, \cite[Theorem 16]{Po08a}
and Lemma \ref{4.1})}.

\smallskip
\item $\mathsf s(C_{3\times 2^a}^3)=8(3\times 2^a-1)+1$, where $a\geq 1$ $($\cite[Theorem
1.8]{G-H-S-T07}$)$.

\smallskip
\item If $G$ is a  $p$-group for some odd prime $p$ with $\mathsf
D(G)=2\exp(G)-1$,
then $\mathsf s(G)=4\exp(G)-3$ $($\cite[Theorem 1.2]{Sc-Zh10a}$)$.

\smallskip
\item If there exists some odd
$q \in \mathbb P$ such that $\mathsf{D}(G_q)-  \exp(G_q)+1 \mid
\exp(G_q)$ and $G_p$ is cyclic for each $p \in \mathbb P \setminus
\{q\}$, then $\mathsf{s}(G) = 2(\mathsf{D}(G_q)-\exp(G_q)) + 2 \exp
(G) - 1$
 $($\cite[Theorem 4.2]{Ge-Gr-Sc11c}; {\rm $G_p$ denotes the $p$-Sylow subgroup of $G$)}.
\end{itemize}
}

\medskip
This shows that precise results for $\mathsf s (G)$ are  extremely
sparse (a few more precise results and upper bounds for groups $G$
which are not of the form $C_n^r$ can be found in
\cite{E-E-G-K-R07,Ge-Gr-Sc11c}). In the present paper, we focus on
groups of the form  $G = C_n^r$,  with $n, r \in \N$ and $n \ge 2$,
and we tackle the study of $\mathsf s (G)$ with a new approach,
combining the direct problem with the associated inverse problem. We
outline this in the next paragraph.

\smallskip
Let $G = C_n^r$ with $n, r \in \mathbb N$ and $n \ge 2$. The inverse
problem associated with $\mathsf s (G)$ asks for the structure of
sequences of length $\mathsf s (G)-1$ that do not have a zero-sum
subsequence of length $n$. The standing conjecture is that every
group of above form satisfies the following Property {\bf D} (see
\cite[Conjecture 7.2]{Ga-Ge06b}).

\begin{propertyd}

 Every sequence \ $S$ \ over $G$ of length \
       $|S| =
      \mathsf s (G) - 1$ \ that has  no zero-sum subsequence of length
      $n$ has the form \ $S = T^{n-1}$ \ for some sequence \ $T$ \ over $G$.
\end{propertyd}

In the case $r=2$, Property {\bf D} was first studied by the second
author in \cite{Ga00a}, and only recently W.A. Schmid completely determined the structure of the sequences having Property {\bf D}
(it was even done for general groups of rank two; see \cite[Theorem
3.1]{Sc11a}). A detailed  overview of Property {\bf D} and its
relationship with further inverse problems can be found in the
survey paper \cite[Section 5]{Ge09a}.

Suppose that  $G = C_n^r$ satisfies Property {\bf D}. Then $\mathsf
s (G) = c (n-1)+1$ where $c = |T|$, and we say that $G$ satisfies
Property {\bf D} with respect to $c$. If  $\mathsf s (G) = c
(n-1)+1$  for some $c \in \N$, then $G$ satisfies the following
Property {\bf D0}.
\begin{propertyd0}(with respect to some $c \in \N$).
Every sequence \ $S$ \ over $G$ of the form $S = g T^{n-1}$ has a zero-sum subsequence of length $n$, where $g \in G$ and $T$ is a sequence of length $|T| = c$.
\end{propertyd0}

Now we can state our main result.

\medskip
\begin{thm}\label{thm1}
Suppose that $C_{m}^{r}$ has Property {\bf D} with respect to $c$
and that $C_n^r$ has Property {\bf D0} with respect to $c$, where
$m, n, r, c \in \N$. If \ $\mathsf s(C_{n}^{r})\leq c(n-1)+n+1$,
\[
n \geq (c-1)^2+1 \quad \text{and} \quad m\geq
\frac{(c(n-1)+n)(n-1)(n^r-(c-1))-(c-1)^2}{n-(c-1)^2} \,,
\]
then
\[
\mathsf s(C_{mn}^{r})\leq c(mn-1)+1 \,.
\]
\end{thm}

\medskip
The proof of Theorem \ref{thm1} will be given in Section \ref{3}.
After the proof we will discuss how to apply Theorem \ref{thm1}, and
we will provide an explicit  list of groups satisfying the
assumptions of Theorem \ref{thm1}. For all of them we will get that
$\mathsf s(C_{mn}^{r}) = c(mn-1)+1$.

\bigskip
\section{Preliminaries} \label{2}
\bigskip

Our notation and terminology are consistent with \cite{Ga-Ge06b} and
\cite{Ge09a}. We briefly gather some key notions and fix the
notation concerning sequences over finite abelian groups.  Let
$\mathbb{N}$ denote the set of positive integers, $\mathbb P \subset
\N$ the set of prime numbers  and
$\mathbb{N}_{0}=\mathbb{N}\cup\{0\}$.
 For real numbers $a, b\in \mathbb R$, we set $[a, b]=\{x \in \mathbb{Z} \mid a\leq x\leq b\}$.
Throughout this article, all abelian groups will be written
additively, and for $n \in \mathbb N$, we denote by $C_n$ a cyclic
group with $n$ elements.

Let $G$ be a finite abelian group and $\exp(G)$ its exponent. A
sequence $S$ over $G$ will be written in the form
\[
S = g_1 \cdot \ldots \cdot g_l = \prod_{g \in G} g^{\mathsf v_g (S)}
\,, \quad \text{with} \ \mathsf v_g (S) \in \mathbb N_0 \ \text{for
all} \ g \in G \,,
\]
and we call
\[
|S| = l \in \mathbb N_0 \quad  \text{the {\it length} \ and } \quad
\sigma (S) = \sum_{i=1}^l g_i = \sum_{g \in G} \mathsf v_g (S)g \in
G \quad  \text{the {\it sum} of} \ S \,.
\]
 The
sequence $S$ is called  a {\it zero-sum sequence} if $\sigma(S)=0$.
For every element $g\in G$, we set
$g+S=(g+g_{1})\cdot\ldots\cdot(g+g_{l})$. Every map of abelian
groups $\varphi \colon G \rightarrow H$ extends to a map from the
sequences over $G$ to the sequences over $H$ by setting
$\varphi(S)=\varphi(g_{1})\cdot\ldots\cdot\varphi(g_{l})$. If
$\varphi$ is a homomorphism, then $\varphi(S)$ is a zero-sum
sequence if and only if $\sigma(S)\in \Ker(\varphi)$.

\medskip
\begin{lem} \label{lem2}
Let $G$ be a finite abelian group.
\begin{enumerate}
\item $\mathsf s(G)\leq |G|+\exp(G)-1$.

\smallskip
\item If $H \subset G$ is a subgroup with $\exp (G) = \exp (H) \exp (G/H)$, then
\[
      \mathsf s (G) \le ( \mathsf s (H) - 1 ) \exp (G/H) + \mathsf s (G/H) \,.
      \]
\end{enumerate}
\end{lem}

\begin{proof}
1. This was first proved by the second author in his thesis (in Chinese). A proof can also be found in \cite[Theorem 4.2.7]{Ge09a}.

\smallskip
2. See \cite[Proposition 5.7.11]{Ge-HK06a}.
\end{proof}

\medskip
\begin{lem} \label{lem4}
Let $n \in \N$ with $n \ge 2$.
\begin{enumerate}
\item $\mathsf s(C_n^r)\geq 2^r(n-1)+1$ for every $r \in \N$.

\smallskip
\item
If $n$ is  odd,  then $\mathsf s(C_{n}^{3})\geq 9n-8$ and $\mathsf s(C_n^4)\geq 20n-19$.
\end{enumerate}
\end{lem}

\begin{proof}
1. See \cite[Hilfssatz 1]{Ha73}.

\smallskip
2. See \cite{El04} and  \cite[Lemma 3.4 and Theorem 1.1]{E-E-G-K-R07}.
\end{proof}

\smallskip
The above mentioned lower bounds for $\mathsf s(C_{n}^{3})$ and
$\mathsf s(C_n^4)$ are due to C. Elsholtz and Y. Edel et al. The
standing conjecture is that equality holds for all odd integers (see
also \cite{G-H-S-T07}).

\medskip
\begin{lem}  \label{lem6}
Let $G=C_{mn}^r$ with $m,n,r \in \mathbb{N}$ and let $c\in
\mathbb{N}$.
\begin{enumerate}
\item  If both $C_m^r$ and $C_n^r$ have Property {\bf D} with respect to $c$ and $\mathsf s (G) = c(mn-1)+1$,
then $G$ has Property {\bf D}.

\smallskip
\item If both $C_m^r$ and
$C_n^r$ have Property {\bf D0} with respect to  $c$, then $G$
 has Property {\bf D0} with respect to  $c$.
\end{enumerate}
\end{lem}

\begin{proof}
1. See \cite[Theorem 3.2]{Ga-Ge-Sc07a}.

\smallskip
2. Let $S=g_{0}\prod_{i=1}^{c}g_{i}^{mn-1}$ be a sequence over
$C_{mn}^{r}$. We  need to show that $S$ has a zero-sum subsequence
of length $mn$.

Let $\varphi \colon G \to G$ denote the multiplication by $m$. Then
$\Ker(\varphi) \cong C_m^r$, $\varphi (G) = mG \cong C_n^r$, and
$$
\varphi(S)=\varphi(g_{0})\prod_{i=1}^{c}\varphi(g_{i})^{mn-1}
$$
is a sequence over $\varphi (G)$. For every  $i \in [1,c]$ and every
$j\in [1,m-1]$, we set $S_{(i-1)(m-1)+j}=g_{i}^{n}$. For the
sequence $T =
S(\prod_{i=1}^c\prod_{j=1}^{m-1}S_{(i-1)(m-1)+j})^{-1}$ we get
$\varphi (T) = \varphi(g_{0}) \prod_{i=1}^{c}\varphi(g_{i})^{n-1}$,
and since $\varphi (G)$ has  Property {\bf D0}, $T$ has a
subsequence $S_0$ such that $\varphi(S_0)$ is a zero-sum sequence of
length $n$. Since $\Ker ( \varphi)$ has Property {\bf D0} and
\[
\prod_{k=0}^{c(m-1)}\sigma(S_k)=\sigma(S_0)\prod_{i=1}^c\prod_{j=1}^{m-1}\sigma(S_{(i-1)(m-1)+j})
=\sigma(S_0)\prod_{i=1}^c( n g_i )^{m-1}
\]
is a sequence over $\Ker (\varphi)$, it
 has a zero-sum subsequence
of length $m$. Therefore there is a subset $I\subset [0,c(m-1)]$
such that $|I|=m$ and $\sum_{k\in I}\sigma(S_k)=0$, which implies
that $\prod_{k\in I}S_k$ is  a zero-sum subsequence of $S$ of length
$mn$.
\end{proof}

\medskip
\begin{lem}\label{lem7}
Let $a,b\in \mathbb{N}_0$.
\begin{enumerate}
\item  $C_{2^a}^r$ has Property
{\bf D} with respect to $2^r$ for every $r\in \mathbb{N}$.

\smallskip
\item $C_{3^a}^4$ has Property {\bf D} with respect to $20$.

\smallskip
\item  $C_{3^{a}5^{b}}^{3}$ has Property {\bf D} with respect to $9$.
\end{enumerate}
\end{lem}

\begin{proof}
1. Obviously,  $C_2^r$ has Property {\bf D}, and Theorem B shows
that Property {\bf D} holds with respect to $2^r$. Using Lemma
\ref{lem6} and Theorem B again, we infer that $C_{2^a}^r$ has
Property {\bf D} with respect to $2^r$.

\smallskip
2. $C_3^r$ has Property {\bf D}   by  \cite[Hilfssatz 3]{Ha73} and
\cite[Lemma 2.3.3]{E-E-G-K-R07}. It follows from Lemma \ref{lem6}
and Theorem B that $C_{3^a}^4$ has Property {\bf D} with respect to
$20$.

\smallskip
3. As mentioned above, $C_{3}^{3}$ has Property {\bf D}, and Theorem
B shows that Property {\bf D} holds with respect to $9$.  It has
been proved in \cite[Theorem 1.9]{G-H-S-T07} that $C_{5}^{3}$ has
Property {\bf D} with respect to $9$. Thus $C_{3^{a}5^{b}}^{3}$ has
Property {\bf D} with respect to $9$ again by  Lemma \ref{lem6} and
Theorem B.
\end{proof}

\medskip
\begin{lem}\label{lem5}
Let $n \in \N$ be an odd integer   which is only divisible by primes
$p \in \{3,5,7,11,13\}$. Then $C_n^3$ has Property {\bf D0} with
respect to $9$.
\end{lem}

\begin{proof}
By Lemma \ref{lem6}, it suffices to show that $C_p^3$  has Property
{\bf D0} with respect to $9$ for all $p \in \{3,5,7,11,13\}$. For $p
\in \{3,5\}$, this follows from Lemma \ref{lem7}. For the other
primes this has been verified  by a computer program written in C
language (the running time was about 0.03, 17 and 31 computer hours,
respectively).
\end{proof}

\bigskip
\section{Proof of Theorem \ref{thm1} and some applications}
\label{3}
\bigskip

\begin{proof}[Proof of Theorem \ref{thm1}]
Let $G = C_{mn}^r$ with $m,n,r \in \N$, and let all assumptions be
as in Theorem \ref{thm1}. Assume to the contrary, there exists a
sequence $S$ over $G$ with $|S|=c(mn-1)+1$ such that $S$ has no
zero-sum subsequence of length $mn$. Let $\varphi \colon G \to G$
denote the multiplication by $m$. Then $\Ker(\varphi) \cong C_m^r$
and $\varphi (G) = mG \cong C_n^r$. We start with a simple
observation which will be used several times in the proof.

\smallskip
\begin{enumerate}
\item[{\bf A1.}\,] Suppose that $S = T_1 \cdot \ldots \cdot T_{c(m-1)}T'$, where $T_1,
\ldots , T_{c(m-1)}, T'$ are sequences over $G$ and, for every $i
\in [1, c(m-1)]$, $\varphi (T_i)$ has sum zero and length $|T_i| =
\exp (\varphi (G)) = n$. Then
\[
\sigma(T_{1}) \cdot \ldots \cdot
\sigma(T_{c(m-1)})=\prod_{i=1}^{c}a_{i}^{m-1} \,,
\] where  $a_{1},
\ldots, a_{c(m-1)} \in \Ker(\varphi)$ are pairwise distinct.
\end{enumerate}

\smallskip

\noindent {\it Proof of \,{\bf A1}}.\, Since $S$ has no zero-sum
subsequence of length $mn$, the sequence $\sigma(T_{1}) \cdot \ldots
\cdot \sigma(T_{c(m-1)})$ has no zero-sum subsequence of length $m$.
Since $\Ker ( \varphi)$  has Property {\bf D}, the assertion
follows. \qed

\smallskip
First we show that $S$ has a product decomposition as in assertion
{\bf A1}. Note that
\[
|\varphi(S)|=c(mn-1)+1=(c(m-1)-1)n+c(n-1)+n+1 \,.
\]
Since $\mathsf s(C_{n}^{r})\leq c(n-1)+n+1$, $S$ allows a product
decomposition $S = T_1 \cdot \ldots \cdot T_{c(m-1)}T'$, where $T_1,
\ldots , T_{c(m-1)}, T'$ are sequences over $G$ and, for every $i
\in [1, c(m-1)]$, $\varphi (T_i)$ has sum zero and length $|T_i| =
\exp (\varphi (G)) = n$ (for details see \cite[Proposition
5.7.10]{Ge-HK06a}).

We set
\[
\varphi(S)=h_{1}^{r_{1}} \cdot \ldots \cdot  h_{t}^{r_{t}} \quad
\text{and} \quad S=S_{1} \cdot \ldots \cdot S_{t} \,,
\]
where $h_{1}, \ldots, h_{t} \in \varphi (G)$ are pairwise distinct,
$r_{1},  \ldots, r_{t}\in \N$, and  $\varphi(S_{i})=h_{i}^{r_{i}}$
for all $i\in [1, t]$. After renumbering if necessary  there exists
an integer $f\in [0, t]$ satisfying

\begin{equation}
\left\{ \begin{aligned}
         r_{i} &\geq (c(n-1)+n)(n-1)\ \ \ \ \ \ , \ \  \text{if} \ i\in [1, f], \\
                  r_{i}&\leq (c(n-1)+n)(n-1)-1,
\ \ \text{otherwise}.
                          \end{aligned} \right.
                          \end{equation}

\medskip
\begin{enumerate}
\item[{\bf A2.}\,] For every $i \in [1, t]$ we have $r_{i}\leq
mn+c(m-1)-m$, and $f \geq c$.
\end{enumerate}

\smallskip

\noindent {\it Proof of \,{\bf A2}}.\, Assume to the contrary, there
exists some $i\in [1, t]$ such that $r_{i}\geq mn+c(m-1)-m+1$. By
the definition of $S_i$, we have $S_{i}=(g+g_{1}) \cdot \ldots \cdot
(g+g_{r_{i}})$ for some $g\in G$ with $\varphi(g)=h_{i}$ and
$g_{j}\in \Ker(\varphi)$ for every $j\in [1, r_{i}]$. Since $\mathsf
s(C_{m}^{r})=c(m-1)+1$ and $r_i\geq m(n-1)+c(m-1)+1$, we can write
$g_{1}\cdot \ldots \cdot g_{r_{i}} = R_0R_{1} \cdot \ldots \cdot
R_{n}$ where $R_j$ is a zero-sum sequence of length $|R_{j}|=m$ for
every $j\in [1,n]$. Then the shifted sequence $g + R_{1}\cdot \ldots
\cdot R_{n}$ is a subsequence of $S_{i}$ such that $|g+R_{1}\cdot
\ldots \cdot R_{n}|=|R_{1}\cdot \ldots \cdot R_{n}|=mn$ and
$\sigma(g+R_{1}\cdot \ldots \cdot R_{n})=mng+ \sum_{j=1}^{n}\sigma
(R_{j})=0$, a contradiction to the assumption that $S$ has no such
zero-sum subsequence.

Combining the upper bounds on $r_i$ with the assumptions that $m\geq
\frac{(c(n-1)+n)(n-1)(n^r-(c-1))-(c-1)^2}{n-(c-1)^2}$ and $n>
(c-1)^2$, we deduce that the $c$-th largest $r_{i}$ is at least
$$\frac{|S|-(c-1)(mn+c(m-1)-m)}{n^r-(c-1)}\geq (c(n-1)+n)(n-1).$$ Thus it follows  that
$f\geq c$. \qed

\medskip
\begin{enumerate}
\item[{\bf A3.}\,] For every $i\in [1,f]$,
$S_{i}=g_{i}^{v_{i}}W_{i}$ for some $g_i\in G$ and $|W_{i}|\leq 1$.
\end{enumerate}

\smallskip

\noindent {\it Proof of \,{\bf A3}}.\, Let $i \in [1, f]$. Since
$|S_i| = r_{i}\geq (c(n-1)+1)(n-1)>2n$, we may choose  an arbitrary
subsequence $L$ of $S_i$ with $|L|=2n$. We set $L = L_1L_2$ where
$|L_1| = |L_2| = n$, and since $\varphi(S_{i})=h_{i}^{r_{i}}$, it
follows  that $\sigma(L_{1}), \sigma(L_{2})\in \Ker(\varphi)$.

Since $|S|=c(mn-1)+1$, $SL^{-1}$ admits a product decomposition
$SL^{-1} = V_0 V_{1} \cdot \ldots \cdot  V_{cm-c-2}$, where
$|V_i|=n$ and $\sigma(V_{i})\in \Ker(\varphi)$ for all $i\in [1,
cm-c-2]$ (we use again \cite[Proposition 5.7.10]{Ge-HK06a}). Now by
{\bf A1}, $\sigma(L_{1})\sigma(L_{2})\sigma(V_{1})\cdot \ldots \cdot
\sigma(V_{cm-c-2})=\prod_{i=1}^{c}a_{i}^{m-1}$, where all $a_{i}\in
\Ker(\varphi)$ are pairwise distinct. After renumbering if necessary
we may assume that $\sigma(V_{1})\cdot \ldots \cdot
\sigma(V_{cm-c-2})=a_1^{k_1}a_2^{k_2}\prod_{i=3}^ca_i^{m-1}$, where
$k_1,k_2\in [m-3,m-1]$ and $k_1+k_2=2m-4$. Therefore,
$\sigma(L_1),\sigma(L_2) \in \{a_1,a_2\}$.

Since $m\geq 4$, $L_1$ is an arbitrary subsequence of $L$ and $L$ an
arbitrary subsequence of $S_i$, we infer that $L$ and therefore
$S_i$ has at most two distinct elements. Therefore there exists some
element $g_i\in G$ which occurs at least
$\frac{r_i}{2}=\frac{(c(n-1)+n)(n-1)}{2}\geq 2(n-1)$ times in $S_i$.
We set $S_i = g_i^{v_i} W_i$ where $v_i=\mathsf v_{g_i}(S_i)$ and
$W_i = a^{|W_i|}$. Assume to the contrary that $|W_i|\geq 2$.  We
set
\[
L_1=g_i^n \mbox{ and } L_2=g_i^{n-2}a^2 \,,
\]
and as above we obtain a product decomposition of $S(L_1L_2)^{-1}$,
say $S(L_1L_2)^{-1} = V_0 V_{1} \cdot  \ldots  \cdot V_{cm-c-2}$,
where $|V_i|=n$ and $\sigma(V_{i})\in \Ker(\varphi)$ for all $i\in
[1, cm-c-2]$. Now by {\bf A1},
$\sigma(L_{1})\sigma(L_{2})\sigma(V_{1})\cdot \ldots \cdot
\sigma(V_{cm-c-2})=\prod_{i=1}^{c}a_{i}^{m-1}$, where all $a_{i}\in
\Ker(\varphi)$ are pairwise distinct. Let $L_1'=L_1ag_i^{-1}$ and
$L_{2}'=L_{2}g_ia^{-1}$.  Since $m\geq 4$, again by {\bf A1}, we
infer that
\[
\sigma(L_{1}')\sigma(L_{2}')\sigma(V_{1})\cdot \ldots \cdot
\sigma(V_{cm-c-2})=\sigma(L_{1})\sigma(L_{2})\sigma(V_{1})\cdot
\ldots \cdot \sigma(V_{cm-c-2})=\prod_{i=1}^{c}a_{i}^{m-1} \,.
\]
 It
follows that $\{\sigma(L_{1}), \sigma(L_{2})\}=\{\sigma(L_{1}'),
\sigma(L_{2}')\}$, which implies that $\sigma(L_1)=\sigma(L_{1}')$
or $\sigma(L_1)=\sigma(L_{2}')$, and thus $g_i=a$, a contradiction.
\qed

\bigskip
Now we have
\[
S = g_{1}^{v_{1}} \cdot \ldots \cdot  g_{f}^{v_{f}} T \quad
\text{where} \quad T=W_{1}\cdot \ldots \cdot W_{f}S_{f+1}\cdot
\ldots \cdot S_{t} \,.
\]

\medskip
\begin{enumerate}
\item[{\bf A4.}\,] $|\supp(\sigma(g_{1}^n)\cdot \ldots \cdot \sigma(g_{f}^n))|\geq
c$.
\end{enumerate}

\smallskip

\noindent {\it Proof of \,{\bf A4}}.\, Assume to the contrary that
$|\supp(\sigma(g_{1}^n)\cdot \ldots \cdot \sigma(g_{f}^n))|\leq
c-1$. By the definition of $f$ we have that $|T|=|W_1\cdot \ldots
\cdot W_f|+|S_{f+1}\cdot \ldots \cdot S_t| \leq
f+[(c(n-1)+n)(n-1)-1](n^r-f)$. Since
$$m\geq \frac{(c(n-1)+n)(n-1)(n^r-(c-1))-(c-1)^2}{n-(c-1)^2}\geq
\frac{n^{r+1}+2n+(c(n-1)+n)(n-1)n^r-cn+c-1}{n},$$  a
straightforward calculation shows that
\[
|T|\leq
(c(mn-1)+1)-[(c-1)(m-1)+1+f]n \,.
\]
Thus we get $v_{1}+\ldots +v_{f}\geq
((c-1)(m-1)+1+f)n$, and hence
\[
\Big\lfloor \frac{v_{1}}{n} \Big\rfloor +\ldots
+ \Big\lfloor \frac{v_{f}}{n} \Big\rfloor \geq (\frac{v_{1}}{n}-1)+\ldots
+(\frac{v_{f}}{n}-1)=\frac{v_{1}+\ldots v_{f}}{n}-f\geq
(c-1)(m-1)+1 \,.
\]
By the pigeonhole principle, there are at least $m$
sequences $C_{1}, \ldots, C_{m}$  among of the
$\lfloor \frac{v_{1}}{n} \rfloor + \ldots +  \lfloor \frac{v_{f}}{n} \rfloor$ sequences
$$\underbrace{g_1^n,
\ldots, g_1^n}_{\lfloor \frac{v_1}{n} \rfloor},\underbrace{g_2^n, \ldots,
g_2^n}_{\lfloor \frac{v_2}{n} \rfloor}, \ldots , \underbrace{g_f^n, \ldots,
g_f^n}_{ \lfloor \frac{v_f}{n} \rfloor}$$
 such that $\sigma(C_{1})=\ldots
=\sigma(C_{m})$. This implies that $C_{1} \cdot \ldots  \cdot C_{m}$ is a zero-sum
subsequence of $S$ of length $mn$, a contradiction. \qed

\bigskip
After renumbering if necessary we may suppose that
$|\supp(\sigma(g_{1}^n) \cdot \ldots \cdot  \sigma(g_{c}^n))|=c$.
Let $Q$ be the subsequence of $S$ with $\varphi(Q)=h_{c+1}^{r_{c+1}}
\cdot \ldots \cdot h_{t}^{r_{t}}$. Then we get
$\varphi(S)=h_{1}^{r_{1}}\cdot \ldots \cdot
h_{c}^{r_{c}}\varphi(Q)$, and we distinguish two cases.

\smallskip
\noindent {\bf Case 1.} $h_{1}^{n-1} \cdot \ldots \cdot h_{c}^{n-1}$
has no zero-sum subsequence of length $n$.

Let $l \in \mathbb N_0$ be maximal such that $Q$ admits a product
decomposition of the form $Q = Q' U_1 \cdot \ldots \cdot U_l$, where
$|U_i|=n$ and $\varphi(U_i)$ is a zero-sum sequence for every $i\in
[1,l]$. It follows that
\[
|Q'|=|\varphi(Q(\prod_{i=1}^{l}U_{i})^{-1})|\leq \mathsf s( \varphi
(G) )-1\leq c(n-1)+n \,.
\]
Since  $\varphi (G) \cong C_n^r$ has Property {\bf D0} with respect
to $c$, every sequence of the form $h_{1}^{n-1} \cdot \ldots \cdot
h_{c}^{n-1}\varphi(x)$ with $ x\in Q'$ has a zero-sum subsequence of
length $n$. Thus for every $x\in \supp (Q')$, one can find a
sequence $U_{l+1} = x U_{l+1}'$, where  $U_{l+1}' \mid SQ^{-1}$,
$|U_{l+1}|=n$ and $\varphi(U_{l+1})$ has sum zero. Since
\[
r_{i}\geq (c(n-1)+n)(n-1)\geq (n-1)|Q'| \quad \text{ for all} \quad
i\in [1, c] \,,
\]
we can do so for every $i\in [1,|Q'|]$, and we obtain a product
decomposition $S = Q'' U_1\cdot \ldots \cdot U_lU_{l+1} \cdot \ldots
\cdot U_{l+|Q'|}$ where the sequences $U_{l+1}, \ldots, U_{l+|Q'|}$
have the above properties. Obviously, we have
$\varphi(Q'')=h_1^{q_1}\cdot \ldots \cdot h_c^{q_c}$. Next we choose
$\lambda= [\frac{q_1}{n}]+\ldots +[\frac{q_c}{n}]$
 subsequences $U_{l+|Q'|+1}, \ldots, U_{l+|Q'|+ \lambda}$ of $Q''$ such that $\varphi (U_{l+|Q'|+i}) \in \{h_1^n, \ldots, h_c^n\}$ for all $i \in [1, \lambda]$
and  $S = Q'''\prod_{i=1}^{l+|Q'|+\lambda}U_i$. Obviously, we have
$\varphi(Q''')=h_1^{q_1'}\cdot \ldots \cdot h_c^{q_c'}$ with
$q_i' \in [0, n-1]$ for all $i\in [1,c]$. Therefore we get
\[
l+|Q'|+\lambda =\frac{|S|-|Q'''|}{n}\geq \frac{c(mn-1)+1-c(n-1)}{n}
\ge  c(m-1)+1 = \mathsf s(C_m^r) \,.
\]
Since $\sigma(U_i)\in \Ker(\varphi) \cong C_m^r$ for all $i\in
[1,l+|Q'|+\lambda]$, the sequence
$\prod_{i=1}^{l+|Q'|+\lambda}\sigma(U_i)$ has a zero-sum subsequence
of length $m$, and hence $S$ has a zero-sum subsequence of length
$mn$, a contradiction.

\smallskip
\noindent {\bf Case 2.} $h_{1}^{n-1}\cdot \ldots \cdot h_{c}^{n-1}$
has a zero-sum subsequence of length $n$.

Let $$h_{1}^{x_{1}}\cdot \ldots \cdot h_{c}^{x_{c}}$$ be a zero-sum
subsequence of $h_{1}^{n-1}\cdot \ldots \cdot h_{c}^{n-1}$ with
$x_{i}\in [0, n-1]$ and $x_{1}+\ldots +x_{c}=n $. Then
$\sigma(g_{1}^{x_{1}}\cdot \ldots \cdot g_{c}^{x_{c}})\in
\Ker(\varphi)$. Since $|\supp(\sigma(g_{1}^{n})\cdot \ldots \cdot
\sigma(g_{c}^{n}))|=c$,  {\bf A1} implies that
$\sigma(g_{1}^{x_{1}}\cdot \ldots \cdot
g_{c}^{x_{c}})=\sigma(g_{k}^n)$ for some $k\in [1, c]$, say $\sigma(g_{1}^{x_{1}}\cdot \ldots \cdot g_{c}^{x_{c}})=\sigma(g_{1}^n)$. Next we  write $S$ in the form
$$S=(g_{1}^{n})^{s_{1}}\cdot \ldots \cdot
(g_{c}^{n})^{s_{c}}g_{1}^{y_{1}}\cdot \ldots \cdot g_{c}^{y_{c}}M,$$
where $s_{i}\in \mathbb{N}$ and $y_{i}\in [0, n-1]$ for all $i\in
[1, c]$. We set $$M_{1}=g_{1}^{y_{1}}\cdot \ldots \cdot g_{c}^{y_{c}}M
\mbox{ and } M_{2}=(g_{1}^{n})^{s_{1}}\cdot \ldots \cdot
(g_{c}^{n})^{s_{c}} \,,$$ and consider $M_{2}$ as a product of $s_1+s_2+\ldots +s_c$
 subsequences of the form $g_1^{n}, \ldots, g_c^{n}$.
 On the other hand, $M_1$ admits a product decomposition of the form
\[
M_1 = M_1' A_1 \cdot  \ldots \cdot A_{c(m-1)-(s_1+s_2+\ldots +s_c)}
\]
such that $|A_i|=n$ and $\sigma(A_i)\in \Ker (\varphi)$ for all $i \in [1, c(m-1)-(s_1+\ldots +s_c)]$.
Since $\sigma(g_{1}^{n}),\ldots, \sigma(g_{c}^{n})$ are pairwise distinct,
{\bf A1} implies that the sequence $\sigma(A_1) \cdot \ldots \cdot
\sigma(A_{c(m-1)-(s_1+\ldots +s_c)})$ contains the element
$\sigma(g_1^n)$ exactly $m-1-s_1$ times. By renumbering if necessary
we assume that
$$\sigma(A_{j})=\sigma(g_{1}^{n})$$  for every  $j\in [1, m-1-s_{1}]$.

Next we provide a further construction of more than $s_1$ subsequences of
$M_{2}$ of length $n$ and with sum $\sigma(g_1^n)$, which
allows us to find more than $s_1$ such  subsequences and
derive a contradiction.
 Since $\sigma(g_{1}^{x_{1}}\cdot \ldots \cdot
g_{c}^{x_{c}})=\sigma(g_1^n)$ and $s_{i}>x_{i}$, we can write $M_2$ in the form $M_2 = M_2'B_1 \cdot \ldots \cdot B_n$ where $B_1 = \ldots = B_n =  g_{1}^{x_{1}}\cdots g_{c}^{x_{c}}$.  Next we write $M_2' = M_2''B_{n+1} \cdot
\ldots \cdot B_{n+[\frac{ns_{1}-nx_{1}}{n}]}$ where
$B_{j}=g_{1}^{n}$ for all $j\in [n+1, n+s_{1}-x_{1}]$.

Thus altogether there are $N=n+s_{1}-x_{1}$  subsequences
$B_{1}, \ldots, B_{N}$ such that $\sigma(B_{j})=\sigma(g_1^n)$ and
$|B_j|=n$ for all $j\in [1, N]$. Since $N=n+s_{1}-x_{1}>s_{1}$, the
sequence $A_1\cdot \ldots \cdot A_{m-s_1-1}B_1\cdot \ldots \cdot
B_{s_1+1}$ is a zero-sum subsequence of $S$ of length $mn$, a
contradiction.
\end{proof}

\medskip
Now we discuss how to apply Theorem \ref{thm1}. Let $r, c$ and $n_0$
be positive integers and $p \in \mathbb P$ a prime. Suppose that
$C_p^r$ has Property {\bf D} with respect to $c$, and  that
$C_{n_0}^r$ has Property {\bf D0} with respect to $c$. By Lemma
\ref{lem6}, $\mathsf s(C_m^r)=c(m-1)+1$ and $C_m^r$ has Property
{\bf D} for
 every $m=p^a$ and every $a\in \mathbb{N}.$  By Lemma \ref{lem2}.2 we get,
 \[
 \mathsf s(C_{mn_0}^{r})\leq n_0(\mathsf s(C_{m}^{r})-1)+\mathsf s(C_{n_0}^{r})=
 n_0c(m-1)+\mathsf s(C_{n_0}^{r})= c(mn_0-1)-n_0+c+\mathsf s(C_{n_0}^{r}) \,.
 \]
 Therefore, for every fixed $n_0$ and $p$, we can choose $a$ sufficiently large such that for $m_0 = p^a$ we get
 $\mathsf s(C_{m_0 n_0}^r)\leq c(m_0 n_0-1)+m_0 n_0+1$ and $m_0 n_0 \ge (c-1)^2 + 1$. Then  we can apply Theorem \ref{thm1}
 with    $n=m_0n_0$ and $m = p^b$ where $b$ is sufficiently large such that $m$ is greater than or equal to the lower bound in $n$.

\smallskip
We work out a few explicit cases.
 Let $a,b,c,d,e$ be nonnegative integers. By the above arguments, we can prove that $\mathsf s(C_{mn}^r)=c(mn-1)+1$ in each of the following situations.

 \begin{enumerate}

 \item Let $r=3, c=9$, $n\geq 65$
an odd  integer such that $C_p^3$ has Property {\bf D0} with
respect to $9$ for all prime divisors $p$ of $n$, and let $m=3^a5^b$
with $$m \geq
\frac{5(n^2-7)\{(50n(n^2-7)-9)(5n(n^2-7)-1)(125n^3(n^2-7)^3-8)-64\}}{(n^2-7)n-64}.$$

\item Let $r=4, c=20$, $n\geq 362$  an odd
integer such that $C_p^4$ has Property {\bf D0} with respect to
$20$ for all prime divisors $p$ of $n$, and let $m=3^a$ with  $$m \geq
\frac{3(n^3-18)\{(63n(n^3-18)-20)(3n(n^3-18)-1)(81n^4(n^3-18)^4-19)-361\}}{(n^3-18)n-361}.$$

\item $r\geq 1, c=2^r$, $n\geq (2^r-1)^2+1$
an even  integer such that  $C_n^r$ has Property {\bf D0}
with respect to $2^r$, and let $m=2^a$ with
$$m\geq
\frac{2n^{r-1}\{(2n^r(2^r+1)-2^r)(2n^r-1)((2n^r)^r-(2^r-1))-(2^r-1)^2\}}{n^r-(2^r-1)^2}$$.

\item Let $r=3, c=9$, $n=7^c11^d13^e \geq 65$, and let $m=3^a5^b$ with  $$m \geq
\frac{5(n^2-7)\{(50n(n^2-7)-9)(5n(n^2-7)-1)(125n^3(n^2-7)^3-8)-64\}}{(n^2-7)n-64}.$$
\end{enumerate}

\begin{proof}
1. By Lemma \ref{lem4}.2, $\mathsf s(C_k^3)\geq 9k-8$ for all
odd positive integers $k$. So, it suffices to prove the upper bound.
Let $a_0,b_0\in \mathbb{N}_0$ with $a_0\in [0,a]$ and $b_0\in [0,b]$
such that,
$$
n^2-7\leq 3^{a_0}5^{b_0}<5(n^2-7).
$$

Let $m_0=3^{a_0}5^{b_0}$ and $n'=m_0n$. By Lemma \ref{lem2} and
Lemma \ref{lem2}.2, $\mathsf s(C_{n'}^3)\leq n(\mathsf
s(C_{m_0}^3)-1)+\mathsf s(C_n^3)=n(9m_0-9)+\mathsf s(C_n^3) \leq
9m_0n-9n+n^3+n-1\leq 9(n'-1)+n'+1$ (the last inequality holds
because $m_0=3^{a_0}5^{b_0}\geq n^2-7$). Let $m'=\frac{m}{m_0}$. By
Lemma \ref{lem7}, $C_{m'}^3$ has Property {\bf D}, and by  Lemma
\ref{lem6}.2
 $C_{n'}^{3}$ has Property {\bf D0} with respect to 9. Now 1. follows from
Theorem \ref{thm1} with $n'$ replacing $n$ and $m'$ replacing $m$.

2. can be proved in a similar way to 1. and we omit it in detail.

3. can be proved in a similar way to 1. by using Lemma \ref{lem4}.1
and we omit it in detail.

4.  It follows from 1. and  Lemma \ref{lem5}.
\end{proof}

\bigskip
\section{Concluding remarks and open problems} \label{4}
\bigskip

We recall the relationship between the Erd{\H o}s-Ginzburg-Ziv
constant and the maximal size of caps in the affine space over
$\mathbb F_3$.

Let $G$ be a finite abelian group, and let $\mathsf g(G)$ denote the
smallest integer $l \in \mathbb N$ such that every squarefree
sequence $S$ over $G$ (or in other terms, every subset $S \subset
G$) of length $|S| \ge l$ has a zero-sum subsequence $T$ of length
$|T| = \exp (G)$. The constant $\mathsf g (G)$ has been studied for
groups of rank two (see \cite{Ga-Th04a} and \cite[Section
5]{Ga-Ge-Sc07a}). Moreover, it found a lot of attention because of
its connection to finite geometry, which we summarize below.

\medskip
\begin{prop} \label{4.1}
Let $G$ be a finite abelian group with $\exp (G) = n \ge 2$.
\begin{enumerate}
\item $\mathsf g (G) \le \mathsf s (G) \le \bigl( \mathsf g (G) -
      1 \bigr)(n-1) + 1$. If $G = C_n^r$, with $n
      \ge 2$ and $r \in \mathbb N$, and
      $\mathsf s (G) = \bigl( \mathsf g (G) - 1 \bigr)(n-1) +
      1$, then $G$ has Property {\bf D}.

\smallskip
\item  Suppose that  $G= \mathbb F_{3}^{r}$. Then the maximal size
of a cap in $G$ equals $\mathsf g (G) - 1$, and we have $\mathsf s
(G) = \big(\mathsf g (G) - 1\big) (3-1) + 1 = 2 \mathsf g (G)  - 1$.

\end{enumerate}
\end{prop}

\begin{proof}
1. The first inequality is clear. For the second statement see
\cite[Lemma 2.3]{E-E-G-K-R07}.

\smallskip
2. This was first observed by H. Harborth (\cite{Ha73}). For a proof
in the present terminology see \cite[Lemma 5.2]{E-E-G-K-R07}.
\end{proof}

\medskip
Let $G = C_n^r$ with $n \ge 3$ odd and $r \in \N$. As already
observed in \cite[Section 5]{E-E-G-K-R07}, in all situations known
so far we have $\mathsf s (G) = \bigl( \mathsf g (C_3^r) - 1
\bigr)(n-1)+ 1$, and we would like to formulate this is a conjecture
(obviously, it implies that $C_n^r$ satisfies Property {\bf D0} with
respect to $\mathsf g (C_3^r) - 1$).

\medskip
\begin{conj} \label{4.2}
For all $n \ge 3$ odd and all $r \in \N$, we have $\mathsf s (C_n^r)
= \bigl( \mathsf g (C_3^r) - 1 \bigr)(n-1)+ 1$.
\end{conj}

\smallskip
Finally we consider groups with even exponent. Let  $n, r$ and $a$
be positive integers.  By Theorem B, Lemma \ref{lem2}.2 and Lemma
\ref{lem2}.1 we obtain that
\[
\mathsf s(C_{2^an}^r)\leq
n(2^r(2^a-1))+n^r+n-1=2^r(2^an)+n^r+n-2^rn-1 \,,
\]
and  by Lemma \ref{lem4}.1 we have
\[
2^r(2^an-1)+1\leq \mathsf s(C_{2^an}^r)\leq
2^r(2^an-1)+1+n^r-2^rn+2^r+n-2 \,.
\]
Therefore there exists an $\alpha \in [0,n^r-2^rn+2^r+n-2]$ such
that
\[
\mathsf s(C_{2^an}^r)=2^r(2^an-1)+1+\alpha \quad \text{for
infinitely many} \quad  a\in \mathbb{N} \,.
\]
We are not aware of any even $n$ such that $\mathsf
s(C_n^r)>2^r(n-1)+1$, and end with the following conjecture.

\medskip
\begin{conj} \label{4.3}
For all $n,r \in \mathbb{N}$ we have
\[
\mathsf s(C_{2^an}^r)=2^r(2^an-1)+1 \quad \text{ for all
sufficiently large} \quad  a\in \mathbb{N} \,.
\]
\end{conj}

\bigskip
{\bf Acknowledgments}. This work was supported by the PCSIRT Project
of the Ministry of Science and Technology, and the National Science
Foundation of China. We would like to thank the referees for all
their comments which helped a lot to improve the presentation of the
paper.

\providecommand{\bysame}{\leavevmode\hbox
to3em{\hrulefill}\thinspace}
\providecommand{\MR}{\relax\ifhmode\unskip\space\fi MR }
\providecommand{\MRhref}[2]{%
  \href{http://www.ams.org/mathscinet-getitem?mr=#1}{#2}
} \providecommand{\href}[2]{#2}

\end{document}